\documentclass[preprint,11pt,graphicx]{article}

\usepackage{amsmath}
\usepackage{amsfonts}
\usepackage{graphicx}
\usepackage{amsthm}
\usepackage{epsfig}

\def\R{{\rm I\kern-3.2pt R}}

\newtheorem{theorem}{Theorem}
\newtheorem{example}{Example}
\newtheorem{lemma}{Lemma}
\newtheorem{corollary}{Corollary}

\newtheorem{proposition}{Proposition}
\newtheorem{assumption}{Assumption}
\newtheorem{definition}{Definition}

\newtheorem{remark}{Remark}

\newcommand {\mm}[1] {\ifmmode{#1}\else{\mbox{\(#1\)}}\fi}

\def\squareforqed{\hbox{\rlap{$\sqcap$}$\sqcup$}}
\def\qed{\ifmmode\else\unskip\quad\fi\squareforqed}

\usepackage{graphics}

\def\ba{\begin{array}}
\def\ea{\end{array}}
\def\beann{\begin{eqnarray*}}
\def\eeann{\end{eqnarray*}}
\def\bea{\begin{eqnarray}}
\def\eea{\end{eqnarray}}
\def\beq{\begin{equation}}
\def\eeq{\end{equation}}

\def\ran{{\rm range}}
\def\los{{\rm loss}}
\def\slos{{\rm scaled\_loss}}

\def\BT{\begin{theorem}}
\def\ET{\end{theorem}}
\def\BL{\begin{lemma}}
\def\EL{\end{lemma}}
\def\BP{\begin{proposition}}
\def\EP{\end{proposition}}
\def\BC{\begin{corollary}}
\def\EC{\end{corollary}}
\def\BD{\begin{definition}}
\def\ED{\end{definition}}
\def\BA{\begin{assumption}}
\def\EA{\end{assumption}}
\def\BR{\begin{remark}}
\def\ER{\end{remark}}
\def\BE{\begin{example}}
\def\EE{\end{example}}

\begin{document}


 \title{{A Robust Robust Optimization Result}
 \thanks{ Supported in part by NSF
through grant DMS-0513337 and ONR
     through grant N00014-08-1-0036.}}

\author{Martina Gancarova and
Michael J.~Todd
\thanks{ School of Operations Research
         and Information Engineering, Cornell University,
        Ithaca, New York 14853, USA}}

\maketitle

\begin{abstract}
We study the loss in objective value when an inaccurate objective
is optimized instead of the true one, and show that ``on average''
this loss is very small, for an arbitrary compact feasible region.
\end{abstract}


\section{Introduction}
\label{intro}
This paper is concerned with the loss in objective value
incurred when an inaccurate objective, because of either
uncertainty or misspecification, is optimized instead of the
true one.

\begin{figure}\label{fig1}
\begin{center}
  \includegraphics[width=6cm]{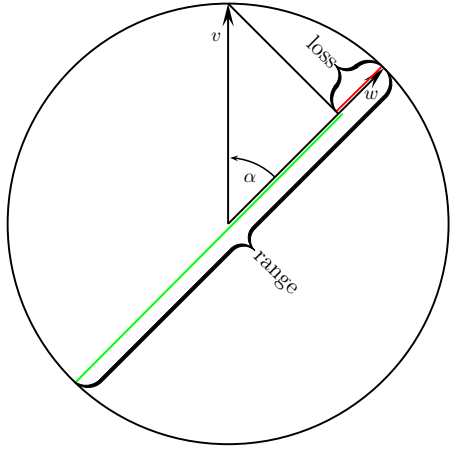}
  \caption{Scaled loss: model case. The loss is the red segment; the range
the sum of the red and green segments.}
\end{center}
\end{figure}

Consider the following model case. Instead of the true objective
$w^Tx$, the nominal objective $v^Tx$ is maximized over the unit
ball in $\R^n$, where $w$ and $v$ are unit vectors making an
angle $\alpha$. (Throughout the paper, we assume
$0 < \alpha < \pi/2$.) Then the computed optimal solution is $x = v$,
attaining a true objective value of $\cos \alpha$. Since the
true optimal value is 1, the loss is $1 - \cos \alpha$, but to
make the measure scale-invariant, we divide by the range of the
true objective over the feasible region, which is 2 (from $-1$
to $+1$). The scaled loss is thus $(1 - \cos \alpha)/2$: see
Figure \ref{fig1}. Our main result claims that this formula for the scaled loss
holds ``on average'' for any compact feasible region. Since this
result on the robustness of the optimal value to misspecification
of the objective holds for any feasible region, we call it a
robust robust optimization result.

Robust optimization has been much studied over the last fifteen
years: see, e.g., \cite{BTEGN1,BTEGN2,BTN1,BTN2,BS}. Usually there is
uncertainty in the constraints as well as the objective, and
the goal is to find a decision vector that is feasible regardless
of the realization of the constraints and that achieves a
guaranteed performance regardless of the realization of the
objective. Typically this leads to an optimization problem that is
harder than the deterministic version of the problem. Our
concerns are appropriate when the decision maker is oblivious to
the error in the objective and does not protect against a
possible misspecification.

In Section 2 we define our setting and give a worst-case bound
on the scaled loss for a class of feasible regions. Section 3
describes two probability distributions for the true and
nominal objectives and obtains our probabilistic result; we
also explain why ``on average'' is in quotes above. Finally,
in Section 4 we discuss the result and outline two applications.

\section{Definitions and Worst-Case Results}

Let $C \subset \R^n$ be compact and nonempty. Since we are
interested in maximizing linear functions over $C$, we could
assume without loss of generality that $C$ is convex, by
replacing it if necessary by its convex hull. However, in Section 4
we treat optimization over a nonlinear transformation of a compact set
and over a set of binary vectors, and so we prefer not to restrict $C$ further.

\BD
For $v \in \R^n$, we define
\beann
\max(v) & := & \max \{v^Tx: x \in C \}, \\
\min(v) & := & \min \{v^Tx: x \in C \}, \mbox{  and} \\
\ran(v) & := & \max(v) - \min(v).
\eeann
\ED

Now consider two objectives, the {\em true} objective $w^Tx$
and the {\em nominal} objective $v^Tx$. If we maximize $v^Tx$
over $C$, the optimal solution set is $\{ x \in C: v^Tx =
\max(v) \}$, and if this is not a singleton, we might be unlucky
and choose the worst $x$ as far as the true objective is
concerned. Hence we make
\BD
For $v, w \in \R^n$, define
\beann
\los(v,w) & := & \max(w) - \min \{w^Tx: x \in C, v^Tx =\max(v) \}
   \mbox{  and} \\
\slos(v,w) & := & \frac{\los(v,w)}{\ran(w)}.
\eeann
\ED
Note that the scaled loss is invariant to translations or
dilations of $C$, and even to rotations if $v$ and $w$ are
correspondingly rotated.

As we have seen, if $C = B^n := \{x \in \R^n: \| x \| \leq 1 \}$
(all norms are Euclidean),
$v = (0;1;0;\dots;0)$, and $w = (\sin \alpha;\cos \alpha;0;\dots;0)$, then $\los(v,w) = 1 - \cos \alpha$ and
$\slos(v,w) = (1 - \cos \alpha)/2$.
On the other hand, if $C$ is the convex hull of $(-1;0;\dots;0)$
and $(+1;0;\dots;0)$ and $v$ and $w$ are as above, then
$\los(v,w) = \ran(w) = 2 \sin \alpha$ and the scaled loss is
1, as bad as it can be.

Note that in this example, the optimal solution set for $v^Tx$ is all of
$C$, and in accordance with the definition above, we choose
the worst of these optimal solutions with respect to the true
objective, namely $(-1;0;\dots;0)$, in evaluating the loss.

Also, for $n = 2$, we will usually view the nominal objective
$v$ as pointing vertically up as in this example. Observe that
there is a subtle bias in this viewpoint. While the decision
maker only sees $v$, and therefore thinks of $v$ as fixed and
$w$ (if she thinks of it at all) as a perturbation of $v$, a
more appropriate perspective would regard
the true objective $w$ as being generated in some suitable
way, and then $v$ arising as a perturbation of $w$.

\begin{figure}[t]\label{fig2}
\begin{center}
  \includegraphics[width=6cm]{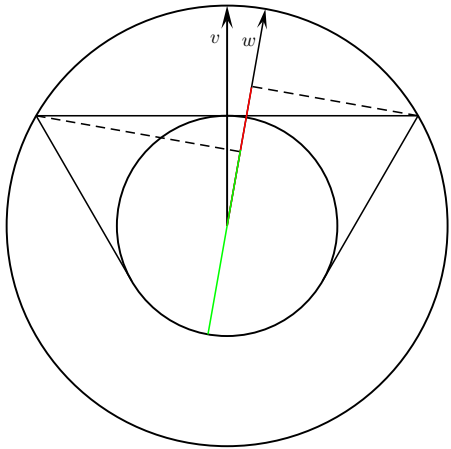}
  \caption{Scaled loss: worst case. The feasible region is the convex hull
of the horizontal line segment and the smaller circle. Again,
the loss is the red segment; the range
the sum of the red and green segments.}
\end{center}
\end{figure}

In the rest of this section, we obtain a worst-case bound on the scaled
loss when $C$ is restricted to avoid the situation
above.
\BT
Assume that $C$ is contained in $B^n$ and contains $rB^n$
for some $0 < r < 1$. Let $v$ and $w$ be two nonzero vectors
making an angle $\alpha$, where $\sin \alpha \leq r \leq \cos \alpha$.
Then, with $\rho := \sqrt{1 - r^2}$, we have
\[
\slos(v,w) \leq \frac{2 \rho \sin \alpha}
   {r(1 + \cos \alpha) + \rho \sin \alpha},
\]
and this bound is tight.
\ET
\begin{proof}
First we show that the right-hand side above can be attained.
Let $n = 2$ and choose $C$ to be the convex hull of
$(-\rho;r)$, $(\rho;r)$, and $rB^2$. Let $v = (0;1)$ and
$w = (\sin \alpha; \cos \alpha)$. Then the set of optimal
solutions for $v^Tx$ is the convex hull of $(-\rho;r)$ and
$(\rho;r)$, with the former being worst for $w^Tx$. With
our assumption that $r \leq \cos \alpha$, the optimal solution
for $w^Tx$ is $(\rho;r)$, while $r \geq \sin \alpha$ implies
that $w^Tx$ is minimized at $-rw$. Hence the loss is
$2 \rho \sin \alpha$ and the range $r(1+\cos \alpha) + \rho
\sin \alpha$, giving the scaled loss as indicated.
See Figure \ref{fig2}.


Now we need to prove the bound. Given any $C$, $v$, and $w$,
we can project $C$ into the plane spanned by $v$ and $w$, and
the projected $C$ will lie in $B^2$ and contain $rB^2$. Hence we
can assume that $n = 2$. By rotating if necessary, we can assume
that $v$ and $w$ are as above (note that the scale of these
vectors is immaterial). Let $s := \max\{v^Tx: x \in C \}
\geq r \geq \sin \alpha$, and let $\sigma := \sqrt{1 - s^2}$.
Then $\min\{w^Tx: x \in C, v^Tx = \max(v)\} \geq - \sigma \sin \alpha
+ s \cos \alpha \geq 0$ and $\min(w) \leq - r$, and so
\[
\slos(v,w) \leq \frac{\max(w) + \sigma \sin \alpha - s \cos \alpha}
  { \max(w) + r};
\]
note that the right-hand side is monotonically increasing in $\max(w)$,
so substituting an upper bound for the latter provides a valid upper
bound on the scaled loss.

If $s \leq \cos \alpha$, then $\max(w) \leq \sigma \sin \alpha +
s \cos \alpha$, and we deduce
\[
\slos(v,w) \leq \frac{2 \sigma \sin \alpha}
  {\sigma \sin \alpha + s \cos \alpha + r} \leq \frac{2 \rho \sin \alpha}
  {\rho \sin \alpha + r \cos \alpha + r},
\]
as desired.

On the other hand, if $s > \cos \alpha$ so that $\sigma < \sin \alpha$,
then $\max(w) \leq 1$ and so
\[
\slos(v,w) \leq \frac{1 + \sigma \sin \alpha - s \cos \alpha}
  {1 + r} \leq \frac{1 + \sin^2 \alpha - \cos^2 \alpha}
  {1 + r} = \frac{2 \sin^2 \alpha}{1 + r} \leq
  \frac{2 \rho \sin \alpha}{r \cos \alpha + \rho \sin \alpha + r},
\]
since $\sin \alpha \leq \rho$ and $1 \geq r \cos \alpha + \rho \sin \alpha$.
Hence the bound is established in either case.

\end{proof}

Since $\rho < \cos \alpha$, the right-hand side above is at most
$\sin \alpha / r$, and it approaches this value for small $r$ and
very small $\alpha$. This bound is of order $\alpha / r$, and hence
much larger than $(1 - \cos \alpha)/2 \approx \alpha^2/4$, the value
in the model case.

\section{Probabilistic Analysis}

Now we examine how the scaled loss behaves when $v$ and $w$ are generated
randomly. We examine two different probability distributions.
We call $v \in \R^n$ a {\em standard Gaussian vector} if its components are
independent standard Gaussian random variables, or equivalently if
$v \sim N(0,I)$.

\BD
We say $(w,v)$ is generated according to Probability Distribution 1
if $w$ and $u$ are independent standard Gaussian vectors in $\R^n$,
and $v := w\cos \alpha + u \sin \alpha$. Expectations with respect
to this distribution are indicated by $E_1$.
\ED

We say a random variable $\xi$ depending on $n$
{\em concentrates} around a positive constant $\beta$ if,
for every positive $\delta$, the probability that $\xi$ lies between
$(1 - \delta) \beta$ and $(1 + \delta) \beta$ converges to 1 as
$n \rightarrow \infty$.
\BP
The angle between $v$ and $w$ generated according to Probability Distribution 1 concentrates around $\alpha$ as $n \rightarrow \infty$.
\EP
\begin{proof}
Because $w$ and $u$ are standard Gaussian random vectors, so is
$v$, and both $w^Tw$ and $v^Tv$ are chi-squared random variables
with $n$ degrees of freedom, both concentrated around their
means, $n$. Also, $v^Tw = \cos \alpha \, w^Tw + \sin \alpha \, u^Tw$,
and since $w^Tw$ is concentrated around $n$ and $u^Tw$ has mean zero
and variance $O(n)$, this is concentrated around $\cos \alpha \, n$.
Hence, using a union bound, we find that with probability approaching 1,
$v^Tw / (v^Tv \, w^Tw)^{1/2}$ lies between $(1 - \epsilon)\cos \alpha /
(1 + \epsilon)$ and $(1 + \epsilon)\cos \alpha / (1 - \epsilon)$ for any
positive $\epsilon$. This implies the result.
\end{proof}

We now define our second model:
\BD
We say $(w,v)$ is generated by Probability Distribution 2 if
$\bar w$ and $\bar u$ are independent standard Gaussian vectors in $\R^n$,
$\hat u = (I - \bar w \bar w^T / \bar w^T \bar w) \bar u$,
$w = \bar w / \| \bar w \|$, $u = \hat u / \| \hat u \|$,
and $v = w \cos \alpha + u \sin \alpha$.  Expectations with respect
to this distribution are indicated by $E_2$.
\ED

Note that all the vectors are well-defined with probability one, with
$w$ and $u$ orthogonal vectors having unit norm, so that
$w$ and $v$ are unit vectors making an angle $\alpha$ with probability one.

In both distributions, $w$ is generated according to some distribution,
and $v$ is generated as a perturbation of $w$. This fits with our interpretation of $w$ as the true objective and $v$ as a nearby inaccurate objective. However, as we now show, we can alternatively
regard $v$ as being generated first and then $w$ as a perturbation
of $v$. This is very useful in our analysis.

\BP\label{pr1}
In both Probability Distribution 1 and
Probability Distribution 2, $(v,w) \sim (w,v)$.
\EP
\begin{proof}
Consider first Probability Distribution 1. Since the matrix
\[
\left[ \ba{cc} \cos \alpha & \sin \alpha \\ \sin \alpha & - \cos \alpha
  \ea \right]
\]
is orthogonal, $(v,z) := (w \cos \alpha + u \sin \alpha,
w \sin \alpha - u \cos \alpha) \sim (w,u)$.
Then, since $(v,w)$ is derived from the first pair
exactly as $(w,v)$ is derived from the second, we see that $(v,w) \sim (w,v)$.

Next assume $(w,v)$ is generated according to Probability Distribution 2.
Then $w$ lies on the unit $(n-1)$-dimensional
sphere $S^{n-1} := \{ x \in \R^n: \|x\| = 1 \}$,
and $v = w \cos \alpha + u \sin \alpha$, where $u$ lies on the unit
$(n-2)$-dimensional sphere
$S_w^{n-2} := \{ x \in \R^n: w^Tx = 0, \|x\| = 1 \}$.

Let $Q$ be any orthogonal matrix. Since
$(Q \bar w,Q \bar u) \sim (\bar w,\bar u)$,
and $Q \bar w / \| Q \bar w \| = Q w$, we see that the distribution
of $Qw$ coincides with that of $w$, so that $w$ is uniformly distributed
on $S^{n-1}$. Also, $(I - Q \bar w \bar w^T Q^T / (Q \bar w)^T(Q \bar w))
Q \bar u = Q (I - \bar w \bar w^T / \bar w^T \bar w) \bar u = Q \hat u$,
so that $(Q \bar w, Q \bar u)$ gives rise to $Q w$ and $Q u$ and hence
$Q v = (Q w) \cos \alpha + (Q u) \sin \alpha)$. It follows that $v$ is
also distributed uniformly on $S^{n-1}$.

Moreover, even if $Q$ is taken to be an orthogonal matrix
that fixes $w$, so that $Q$ depends on $w$, we still have $Q \bar w =
\bar w$ a standard Gaussian vector in $\R^n$ and $Q \bar u$ a standard
Gaussian vector independent of $w$. Hence, proceeding as above, we see
that $(w = Q \bar w, Q \bar u)$ gives rise to $Q \hat u$ and hence $Q u$,
so since $(\bar w, Q \bar u) \sim (\bar w, \bar u)$, we find $Qu \sim u$,
from which, conditioned on $w$, $u$ is uniformly distributed on $S_w^{n-2}$.

Finally, we show that $(v,z) :=
(w \cos \alpha + u \sin \alpha, w \sin \alpha - u \cos \alpha) \sim (w,u)$. Since $w = v \cos \alpha + z \sin \alpha$, this will show
that $(v,w) \sim (w,v)$ as desired. We have already shown that
$v$, the first member of the pair, is distributed uniformly on
$S^{n-1}$. We now consider all pairs $(\bar w, \bar u)$ that
give rise to a given $v$. If $Q$ is an orthogonal matrix that fixes
$v$, then $(Q \bar w, Q \bar u)$ also gives rise to the same $v$.
Since $Q$ is orthogonal, the distribution of $(\bar w, \bar u)$,
conditional on this fixed $v$, is invariant under pre-multiplication
of each vector by such a $Q$. As we have seen, under this transformation
$w$ is transformed to $Qw$ and $u$ to $Qu$, and hence
$z := w \sin \alpha - u \cos \alpha$ is transformed to $Qz$. It follows
that $z$, which has unit norm and is orthogonal to $v$, is uniformly
distributed on $S_v^{n-2}$. This concludes the proof.
\end{proof}

We are now ready to analyze the behavior of the scaled loss ``on average''
for our two models. First we investigate the range function:
\BL\label{lem1}
For $i = 1,2$, we have
\[
E_i \max(v) = E_i \max(w), \quad E_i \, \ran(w) = 2 E_i \max(w).
\]
\EL
\begin{proof}
The first equation follows from Proposition \ref{pr1} above, since
$v$ and $w$ have the same distribution. For the second equation, note
that
\[
E_i \, \ran(w) = E_i \max(w) - E_i \min(w) = E_i \max(w) + E_i \max(-w),
\]
and that $E_i \max(-w) = E_i \max(w)$ since under both probabilistic models,
$w$ has a symmetric distribution.
\end{proof}

Next we examine the loss:
\BL\label{lem2}
For $i = 1,2$, we have
\[
E_i \, \los(v,w) = (1 - \cos \alpha) E_i \max(w).
\]
\EL
\begin{proof}
First note that, since $C$ is compact, the convex function
$\max(v)$ is finite everywhere, and hence is differentiable almost
everywhere, with respect to Lebesgue measure and hence with respect
to Probability Distribution 1. But $\max(v)$ is differentiable at
$v$ exactly when the maximum of $v^Tx$ over $C$ is attained at a
single $x$, which we denote by $x_v$. Since this property is invariant under positive scalings of $v$, we see that it holds also for almost all $v$
under Probability Distribution 2 also.

Hence with probability one,
\[
\los(v,w) = \max(w) - w^Tx_v = \max(w) - (v \cos \alpha + z \sin \alpha)^Tx_v,
\]
where, as in the proof of Proposition \ref{pr1}, we let $z := w \sin \alpha
- u \cos \alpha$. Now in either model, $z$ has a symmetric distribution
conditional on $v$, and since with probability one $x_v$ depends only
on $v$, $z^Tx_v$ has mean zero. Hence
\[
E_i \, \los(v,w) = E_i \max(w) - \cos \alpha \, E_i \max(v) - 0,
\]
and the result follows from Lemma \ref{lem1}.
\end{proof}
From Lemmas \ref{lem1} and \ref{lem2}, we immediately deduce
\BT\label{thm1}
For $i = 1,2$, we have
\[
\frac{E_i \, \los(v,w)}{E_i \, \ran(w)} = \frac{1 - \cos \alpha}{2}.
\]
\ET
\qed

Note that we do not have a result on the expected scaled loss, which would
be an expectation of the ratio of the loss to the range, but only on the
ratio of the expectations, which is why we have put ``on average'' in
quotes above.

\section{Discussion and applications}

It seems at first that the theorem of the last section would hold under
much weaker probabilistic assumptions, merely requiring that $w$ and
$u$ have symmetric distributions. Unfortunately, there are two problems
with this. First, what we really need is that $v$ and $z$ have symmetric
distributions, but putting restrictions on $v$ and $z$ conflicts with the
natural interpretation that the true objective $w$ should be generated
first, and then $v$ as a perturbation of $w$. Second, it is crucial that
$\max(v)$ and $\max(w)$ have the same expectation, and this appears hard
to ensure under weaker assumptions: the fact that $(v,w)$ and $(w,v)$
have the same distribution under our two models is key in our
development.

One way in which the result can be generalized is in allowing a random
choice of $\alpha$. Our two models yield vectors $v$ and $w$ making an
angle that either concentrates around $\alpha$ or is exactly $\alpha$.
Instead, we can consider probability distributions
on the triple $(\alpha,w,v)$ as follows: first
$\alpha$ is generated according to an arbitrary distribution
supported on $(0,\pi/2)$; then, conditional on $\alpha$, $w$ and $v$ are
generated according to Probability Distribution 1 or 2. It is easy
to see that all our arguments can be extended by first conditioning
on $\alpha$, and the expected loss divided by the expected range will be
$(1 - E(\cos \alpha))/2$, where the expectation is taken with respect to
the distribution on $\alpha$.

Another generalization allows very general distributions, but changes
the way the objective vector is perturbed. Let $f_j$ be a symmetric probability
density on $\R$ for $j = 1,\dots,n$. For each $j$, draw $w_j$ and $u_j$
independently from $f_j$, and then let $v_j$ be $w_j$ with probability
$\cos \alpha$ and $u_j$ with probability $1 - \cos \alpha$, with all the
choices independent. Let $t \in \R^n$ be defined by $t_j = +1$ if
$v_j = w_j$, $t_j = -1$ if $v_j = u_j$, so each $t_j$ is $+1$ with
probability $\cos \alpha$ and $-1$ with probability $1 - \cos \alpha$.
Then it is clear how $(w,v)$ arises as a function of the triple $(w,u,t)$.
But if we define $z := w + u - v$, then in each component, $z$ agrees
with $w$ ($u$) exactly when $v$ agrees with $u$ ($w$). It follows
that $(v,z,t)$ has the same distribution as $(w,u,t)$, and $(v,w)$
arises from $(v,z,t)$ as does $(w,v)$ from $(w,u,t)$. Hence $(v,w)$
and $(w,v)$ have the same distribution. Moreover, the arguments of the
previous section can be duplicated, and again lead to the result
that the expected loss divided by the expected range is exactly
$(1 - \cos \alpha)/2$.  Under mild conditions on the $f_j$'s,
the angle between $v$ and $w$ concentrates around $\alpha$.  
Note that, in this model, a small fraction of
the components are changed a possibly large amount, while in the previous 
models, each component is changed a small amount.

\begin{figure}[ht]
\begin{minipage}[b]{0.5\linewidth}
\centering
  \includegraphics[width=6cm]{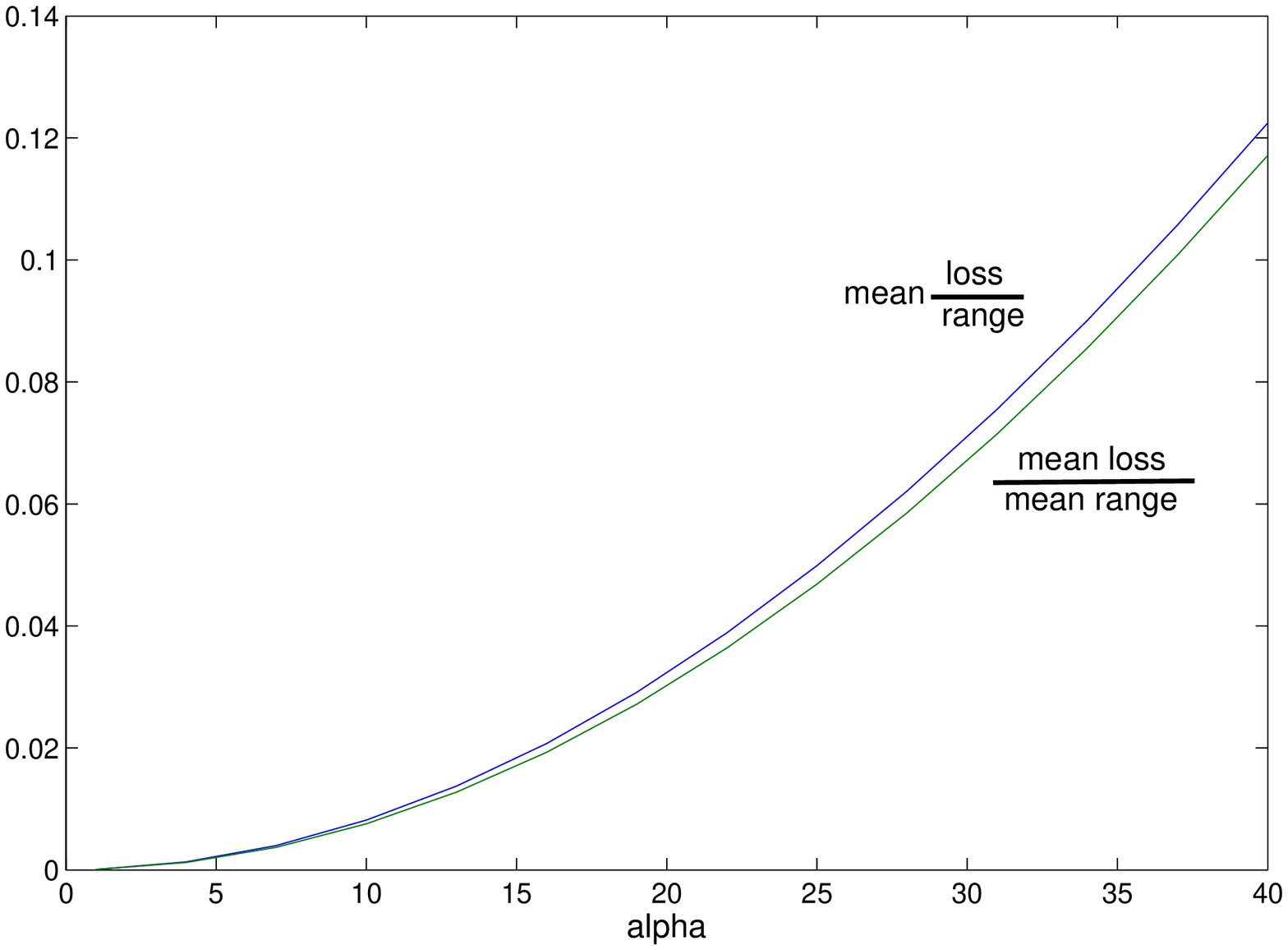}
  \caption{Ratio of expectations versus expectation of ratios for AGG.}
\label{fig3}
\end{minipage}
\hspace{0.5cm}
\begin{minipage}[b]{0.5\linewidth}
\centering
  \includegraphics[width=6cm]{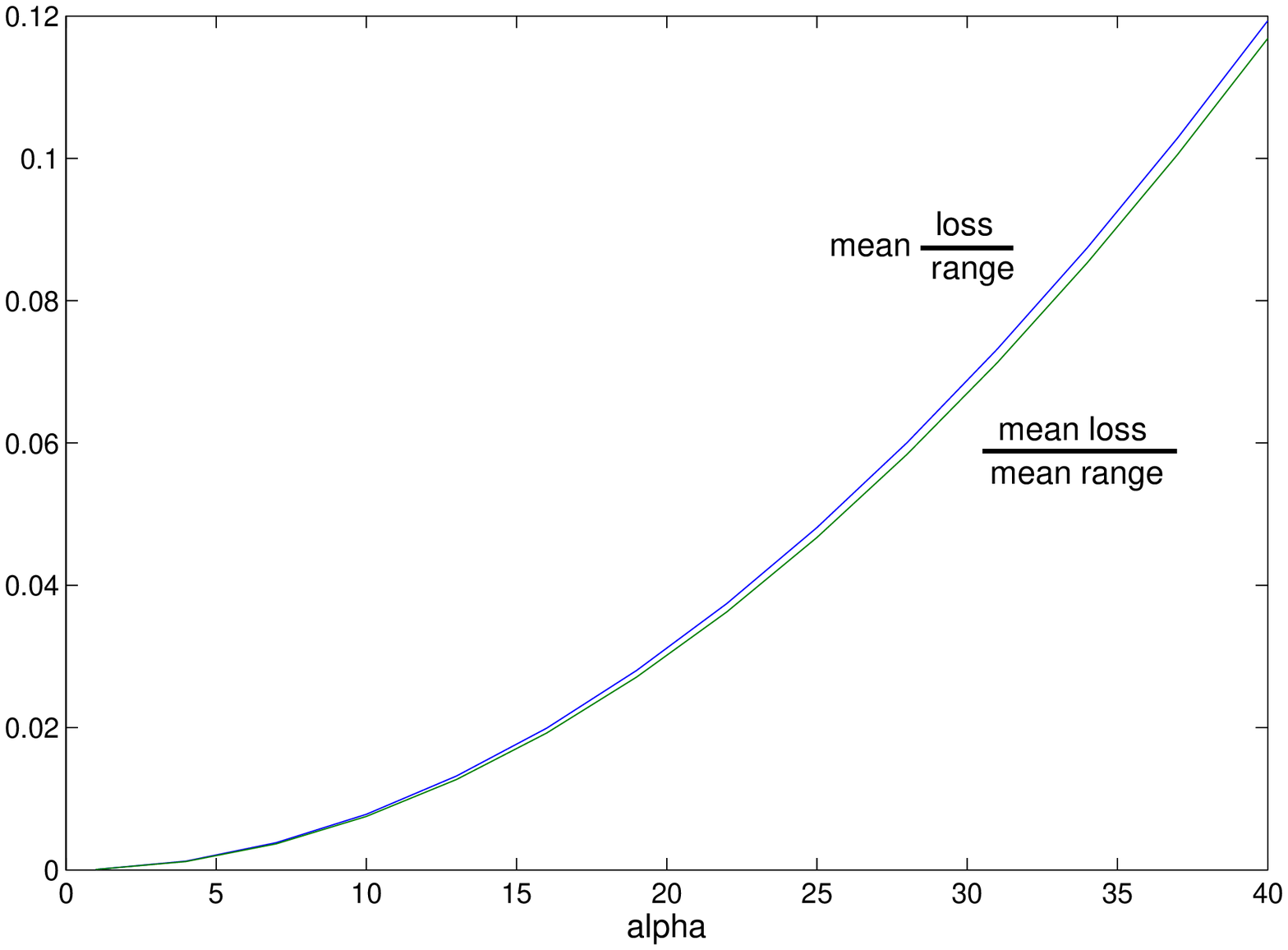}
  \caption{Ratio of expectations versus expectation of ratios for BOEING1.}
  \label{fig4}
\end{minipage}
\end{figure}

We argued in the introduction that the scaled loss provided a good measure
of how much is lost in objective value when implementing the optimal solution
for a misspecified objective. However, the result of the previous section
is concerned with the ratio of the expectations of the loss and the range,
rather than the more meaningful expectation of the ratio. We therefore conducted
some experiments with two NETLIB \cite{Netlib} problems, AGG and BOEING1,
to see how much the results differ. The first has 489 rows and 163 columns,
the second 351 rows and 384 columns. Figures \ref{fig3} and
\ref{fig4}, which give graphs
of the expectation of the ratio and of the ratio of the expectations
as functions of the angle $\alpha$ in degrees, show
that our results should be applicable to the more meaningful measure also.
In both figures, each data point is obtained using at least 10,000 pairs
$(v,w)$ generated from Probability Distribution 1.

%

Our result is limited to clarifying what happens when linear objective functions
are perturbed. In general, little can be said for nonlinear functions,
partially because it is not clear how random nonlinear objective functions
and their perturbations should be defined. However, our analysis can be
applied to one case where objective functions are nonlinear. Suppose there
are several continuous objective functions,
$f_i(x)$ for $i = 1,\dots,k$. The decision
maker is interested in high values of all of these objectives, so we are
in the realm of multi-criteria optimization; see, for instance,
\cite{DasDen,Hill}. Often, a linear combination
$\sum_i w_i f_i(x)$ of the objectives is maximized, and indeed,
any such optimal solution for a positive $w$ is
an efficient or Pareto-optimal solution. (The converse is true when
all $f_i$'s are concave, and in this case the optimization problem is convex,
but this restriction is not needed for our discussion.) However, since the
different function $f_i$ may be hard to compare, it is difficult to decide on
appropriate weights $w$. Our theorem indicates in some sense that the
choice may not matter too much.

Let $Y := \{y = (f_1(x);\dots,f_k(x)): x \in X\}$. Since $X$ is nonempty and
compact, so is $Y$, and trivially $\max \{\sum_i w_i f_i(x): x \in X\}
= \max \{w^Ty: y \in Y\}$. Theorem \ref{thm1} shows that
the latter is insensitive in some precise sense to the specification of
the objective $w$, and this translates directly into the insensitivity
of the original problem's optimal value to the specification of the weights.
Of course, there is a large caveat here: the result requires $v$ and $w$
to be randomly chosen from symmetric distributions, and hence their
components are as likely to be negative as positive, while in the
multi-criteria setting, the weights are always positive. Nevertheless,
we believe our theorem gives some credence to the hope that incorrect
choices of weights should not hurt much.

Our second application is to the complexity of combinatorial optimization
problems. Beier and V\"{o}cking \cite{BV} and R\"{o}glin and Teng
\cite{RT} have uncovered a fascinating connection between binary
optimization problems that can be solved in randomized pseudo-polynomial time,
that is, in randomized polynomial time if the data are encoded in unary, and
smoothed complexity. In particular, R\"{o}glin and Teng show that
such a problem can be solved in expected time polynomial in the input size 
and $1/\sigma$, where the adversarily chosen
objective function coefficients are perturbed
by independent Gaussians with mean zero and variance $\sigma^2$ (there are
slight technical subleties; see Sections 2 and 6 of \cite{RT}).
This is normally interpreted as saying that arbitrarily close to any
potentially hard instance there are polynomially solvable instances.
This gives support to the belief that one would be unlucky to choose
a bad instance, in a rather precise and strong sense.
Our result provides another avenue to solving such problems.
One can {\em explicitly} make a small random
perturbation of the objective function coefficients, thereby obtaining
a problem with provably expected polynomial-time complexity.
(Note that $w + z$, where $z$ has independent zero-mean Gaussian
components with standard deviation $\sigma$, is proportional to
$w \cos \alpha + u \sin \alpha$, where $u$ is a Gaussian random vector
and $\alpha := \arctan \sigma$.) Solving this perturbed
problem gives a feasible solution to the original
problem, and Theorem \ref{thm1}
gives credence to the hope that this solution will be close to optimal
for the true objective function.
(For some models of generating $w$, renormalizing might give a
distribution on $\alpha$, rather than a fixed value, but the extension
mentioned at the beginning of this section allows for this possibility.)
Of course, our result only proves this ``on average," so one would be
unlucky to have an objective function where the loss is large, in a certain
sense. We believe this viewpoint provides further insight into the
notion of smoothed complexity, at least when only the objective function
is perturbed.

\end{document}